\documentclass[11pt,reqno]{amsart}
\usepackage[%
pdftitle={On the spectrum of the absorption operator with bounce-back boundary conditions},
pdfauthor={Khalid Latrach and Bertrand Lods},
pdfsubject={On the spectrum of the absorption operator with bounce-back boundary conditions},
  colorlinks=true,
  linkcolor=blue,
  citecolor=blue,
]{hyperref}
\usepackage[T1]{fontenc}
\usepackage{mathrsfs}
\usepackage{amsmath,amsthm,amssymb}

\usepackage{latexsym}
\usepackage{times}
\usepackage{enumerate}
\usepackage{mathrsfs}
\usepackage{stmaryrd}
\usepackage{amsopn}
\usepackage{amsmath}
\usepackage{amssymb}
\usepackage{amsfonts}
\usepackage{amsbsy}
\usepackage{amscd,indentfirst}
\usepackage{hyperref}
\usepackage{amsfonts,amsmath,latexsym,amssymb,verbatim,amsbsy}
\usepackage{amsthm}
\usepackage{colordvi}
\newtheorem{theo}{\rm \bf Theorem}[section]
\newtheorem{lemme}[theo]{Lemma}
\newtheorem{propo}[theo]{Proposition}

\newtheorem{hyp}[theo]{Assumption}
\newtheorem{defi}[theo]{Definition}

\newtheorem{nb}[theo]{Remark}

\font\teneuf=eufm10 at 12pt \font\seveneuf=eufm7 at 8pt
\font\fiveeuf=eufm5 at 6pt
\newfam\euffam
\textfont\euffam=\teneuf \scriptfont\euffam=\seveneuf
\scriptscriptfont\euffam=\fiveeuf

\newfont{\secgoth}{eufm10 at 16pt}

\def \D {\mathscr{D}}
\def \Omega {\mathcal{D}}
\def \omm {{\Omega} \times \R^N}
\def \ut {(U(t))_{t \geq 0}}

\def \essu {\operatornamewithlimits{ess\,-sup}}
\def \esin {\operatornamewithlimits{ess\,-inf}}

\def \ds {\displaystyle}
\def \vt {(V(t))_{t \geq 0}}

\def \d {\mathrm{d}}
\def\ind#1{\lower5pt\hbox{$\scriptstyle #1$}}

\renewcommand{\epsilon}{\varepsilon}
\renewcommand{\geq}{\geqslant}
\renewcommand{\leq}{\leqslant}

\newcommand{\R}{\mathbb R}
\newcommand{\C}{\mathbb C}

\newcommand{\eps}{\varepsilon}

\newcommand{\n}{\noindent}

\newcommand{\vs}{\vskip 0.3 cm}

\renewcommand{\eps}{\varepsilon}

\renewcommand{\n}{\noindent}

\renewcommand{\vs}{\vskip 0.3 cm}

\numberwithin{equation}{section}
\title[Spectral analysis of transport equations with
bounce-back boundary conditions] {Spectral analysis of transport
equations with bounce-back boundary conditions.}

\author[K. Latrach \& B.  Lods]{}
\begin{document}
\bibliographystyle{plain}
\maketitle \centerline{\scshape K. Latrach \& B. Lods}
\smallskip { \centerline{Universit\'{e} Blaise Pascal (Clermont II)}
\centerline{Laboratoire de Math\'ematiques, CNRS UMR 6620}
\centerline{63117 Aubi\`{e}re,
France.\footnote{\texttt{Khalid.Latrach,Bertrand.Lods@math.univ-bpclermont.fr}}}

\begin{abstract}  We  investigate
the spectral properties  of the time-dependent
linear transport equation with bounce-back boundary conditions. A
 fine analysis  of the  spectrum
of the streaming operator is given and  the explicit expression of
the strongly continuous streaming semigroup is derived.   Next,
making use of a recent result from \cite{Sbihi}, we prove, via a
compactness argument, that the essential spectrum of the transport
semigroup and that  of the streaming semigroup coincide on all
$L^p$-spaces with $1<p<\infty$.

\noindent  \textsc{Keywords:} Transport operator, bounce-back
boudary conditions, transport semigroup, essential spectrum,
compactness.

 \noindent \textsc{AMS subject classifications (2000):}
47D06, 47D05, 47N55, 35F05, 82C40
\end{abstract}

\maketitle

\section{Introduction}

The spectral theory of transport equations with \textit{no-reentry
boundary conditions} (i.e. with zero incoming flux in the spatial
domain) received a lot attention in the last decades (see, for
example, the works \cite{Greiner, Jorgens,Mokhtar1,Mokhtar2,
Takac,Vidav1,Vidav2,Voigt2,Voigt3,Weis} and the references therein).
The picture is fairly complete by now and almost optimal results
have been obtained in \cite{Mokhtar3} for bounded spatial domains
and in \cite{sbihi2} for unbounded domains.

When dealing with reentry boundary conditions (including periodic
boundary conditions, specular reflections,  diffuse reflections,
generalized or mixed type boundary conditions), many progress have
been made in the recent years in the understanding of the spectral
features of one-dimensional models
\cite{Chabi,Dehici,LL,LLM,LM,LS,Lods2}. However, to our knowledge,
for higher dimensions, only few partial results are available in the
literature \cite{Chen,Frosali,Makin,Pal}, dealing in particular with
very peculiar shapes of the spatial domain. Our paper deals with the
following two problems concerning multidimensional transport
equations with \textit{bounce-back (reverse) boundary conditions }in
convex bounded domains:

\begin{enumerate} \item The spectral analysis of the streaming operator subjected to
bounce-back boundary conditions and the explicit expression of the
streaming semigroup.

\item  The compactness of the difference of the (perturbed) transport
  semigroup and the streaming semigroup.
\end{enumerate}

 To be more precise,  we are  concerned with the following
initial-boundary-value problem in $L^p$--spaces $(1 \leq p <
\infty)$
\begin{subequations}\label{1}
\begin{equation}\label{1a}\begin{split}
\displaystyle \frac{\partial \psi}{\partial t}(x,v,t)&=- v
\cdot\nabla_x \psi(x,v,t)-\Sigma (x,v)\psi (x,v,t)+
 \int_{\R^N}\kappa(x,v,v')\psi(x,v',t) dv'\\
 &=T\psi(x,v,t)+K\psi(x,v,t),\qquad \qquad
(x,v)\in \mathcal{D} \times \R^N, \, t> 0;\end{split}\end{equation}
with bounce-back boundary conditions:
\begin{equation}\label{1b}
\psi_{|\Gamma_-}(x, v,t) =\gamma  \psi_{|\Gamma_+}(x,-v,t), \qquad
 (x,v) \in \Gamma_-, t >0;
\end{equation}
and the initial condition
\begin{equation}\label{1c}\psi(x,v,0) =\psi_0(x,v) \in L^p(\mathcal{D} \times
\R^N).\end{equation}\end{subequations} Here $\mathcal{D}$ is a
smooth \textit{convex} open subset of $\R^N$ $(N\geq 1)$, $\gamma$
is a real constant belonging to $(0,1)$  and $\Gamma_\mp$ represent
the incoming and outgoing parts of the boundary of the phase space
(see Section 2 for details). The collision frequency $\Sigma
(\cdot,\cdot)\in L^\infty(\mathcal{D}\times \R^N)$ is a non-negative
function. The scattering kernel $\kappa(\cdot,\cdot,\cdot)$ is
nonnegative and defines the linear operator $K$ called the collision
operator which is assumed to be bounded on $L^p(\mathcal{D}\times
\R^N)$ $(1\leq p < +\infty)$. The operator $T$ appearing in
\eqref{1a} is called the \textit{streaming operator}
 while $T+K$ denotes the (full) transport operator. It is well-known \cite{Beal} that,
 since $\gamma < 1$, $T$ generates a C$_0$-semigroup of contractions $\ut$ in $L^p(\Omega \times \R^N)$ (\textit{streaming semigroup}) and,
since $K$ is bounded, $T+K$ is also the generator of a
C$_0$-semigroup $\vt$ in $L^p(\Omega \times \R^N)$
(\textit{transport semigroup}).\par\vs

Our work is displayed into two parts, referring to the above points $(1)$ and $(2)$.
First, we present a fine description of the spectrum  $\sigma(T)$ of
the streaming operator in $L^p(\Omega \times \R^N)$, $1 \leq p <
\infty.$ Moreover, we derive the explicit expression of the streaming
semigroup $\ut$ for the particular case of a space-homogenous
collision frequency. Second, we prove that the difference
$R_1(t)=V(t)-U(t)$ is compact in $L^p(\Omega \times \R^N)$ $(1 < p < \infty)$ for any
$t \geq 0$ under natural assumptions on the scattering operator $K$.  The
interest of such a compactness result lies in the fact
that it implies that the streaming semigroup and the  transport
semigroup possess the same essential spectrum (see \cite{Schechter} for a
precise definition). In particular,
their essential types coincide. This shows that  the part of the
spectrum of the transport semigroup outside the spectral disc of the
streaming semigroup consists of, at most, eigenvalues with finite
algebraic mutiplicities. Assuming the existence of such eigenvalues,
the transport semigroup can be decomposed into two parts: the first
containing the time development of finitely many eigenmodes, the
second being of faster decay.

Although the well-posedness of the problem (1.1) is a known
 fact \cite{Arlotti,Beal,Frosali,Lods1},  the  description
of the  spectrum of the streaming operator and the analytic
expression of its semigroup seem to be  new. Let us also notice that
besides the interesting consequences of the compactness of $R_1(t)$
on the behavior for large times of the solution of the problem
(1.1),  it is an interesting result in itself. \textit{Actually,  it
is the first time that the compactness of the first order remainder
term of the Dyson-Phillips expansion of the transport operator with
reentry boundary conditions is discussed in higher dimensions}. For
the one-dimensional case we refer to the work \cite{LS} while the
compactness of $R_1(t)$ in the case of non-reentry boundary
condition was established in \cite{Mokhtar3}.\par\vs

As in \cite{LS}, the mathematical analysis is based upon a recent
result owing to M. Sbihi \cite{Sbihi} (see also Section 4) valid for
Hilbert spaces. Actually, under some natural assumptions  on the
collision operator, we prove, via approximation arguments, the
compactness of $R_1(t)$ on $L^2(\mathcal{D} \times \R^N, dx\otimes
dv)$. The result is then extended to $L^p(\mathcal{D} \times \R^N,
dx\otimes dv)$ with $p\in (1,2)\cup (2,\infty)$  by an interpolation
argument. Unfortunately, the limiting case $p=1$ is not covered by
our analysis and requires certainly another approach. Notice that
the compactness of $R_1(t)$ for $p=1$  is also an open problem in
the one-dimensional case \cite{LS} and for no-reentry boundary
condition \cite{Mokhtar3}.
\par\vs

 The outline of this work is as follows. In Section 2 we introduce the functional setting
 of the problem and fix the different notations and facts needed in the sequel.
 Section 3 is devoted to the spectral analysis of the streaming operator with bounce-back
boundary conditions and to the analytic expression of the streaming
semigroup. The compactness  of the first order remainder term of the
Dyson-Phillips expansion is the topic of Section 4.

\section{Preliminary results}

For the definitions of the different spectral notions used
throughout this paper we refer, for example, to the book
\cite{Schechter}. If $X$ is a Banach space, $\mathcal{B}(X)$ will
denote the set of all bounded linear operators on $X$.

 Let $\Omega$ be a
smooth bounded open subset of $\mathbb{R}^N $.  We define the
partial Sobolev space
$$W_p=\{\psi \in X_p\,;\,v \cdot \nabla_x \psi \in X_p\}$$
where $X_p=L^p(\Omega \times \mathbb{R}^N,\d x \otimes \d v)$ $(1
\leq p < \infty)$. Let us denote by $\Gamma_-$ (respectively
$\Gamma_+$) the incoming (resp. outgoing) part of the boundary of
the phase space $\Omega \times \mathbb{R}^N$
$$\Gamma_{\pm}=\left\{(x,v) \;\in \partial \Omega \times \mathbb{R}^N\;;\;\pm v \cdot
n(x) \geq 0\right\}$$ where $n(x)$ stands for the outward normal
unit at $x \in \partial \Omega$. Suitable $L^p$-spaces for the
traces on $\Gamma_{\pm}$ are defined as
$$L^p_{\pm}=L^p\big(\Gamma_{\pm};|v \cdot n(x)|\d\gamma(x)\otimes\d v\big),$$
$\d\gamma(\cdot)$ being the Lebesgue measure on $\partial \Omega$.
For any $\psi \in W_p$, one can define the traces $\psi_{|\Gamma_{\pm}}$
on $\Gamma_{\pm}$, however these traces do not belong to $L^p_{\pm}$ but to
a certain weighted space (see \cite{Beal,Cessenat1,Cessenat2}). For this reason,
one defines
$$\widetilde{W}_p=\bigg\{\psi \in W_p\,;\,\psi_{|\Gamma_{\pm}} \in
L^p_{\pm}\bigg\}.$$  In all the sequel, we shall assume that
$\Sigma(\cdot,\cdot)$ is a measurable \textit{{non-negative}}
function on $\Omega \times \mathbb{R}^N$ that fulfills the
following.

\begin{hyp}\label{asssig} The collision frequency $\Sigma(\cdot,\cdot)$
is an even function of the
velocity, i.e. for any $(x,v) \in \Omega \times \mathbb{R}^N$,
$\Sigma(x,v)=\Sigma(x,-v)$.
\end{hyp}

\n Let us define the absorption operator with \textit{{bounce-back
boundary conditions:}}
\begin{equation*}
\begin{cases}
T\::\: \D(T) &\subset X_p
\longrightarrow X_p\\
&\varphi \longmapsto T\varphi(x,v):=-v \cdot \nabla_x
\varphi(x,v)-\Sigma(x,v)\varphi(x,v),
\end{cases}
\end{equation*}
with domain $$\D(T):=\bigg\{\psi \in \widetilde{W}_p \;\text{ such
that } \psi_{|\Gamma_-}(x,v)=\gamma\, \psi_{|\Gamma_+}(x,
-v)\bigg\}$$ where $0 < \gamma  <1.$ We recall (see e.g.
\cite{Beal,Frosali})  that $T$ is a generator of a  non-negative
C$_0$-semigroup of contractions $\ut$ in $X_p$ $(1 \leq p < \infty)$
\begin{nb} Notice that our analysis also applies to the more general
case $\gamma \geq 1$ provided the associated transport operator $T$
 generates a C$_0$-semigroup in $X_p$. Practical conditions
on the geometry of $\Omega$  ensuring the latter to hold are given
in \cite{Lods1}.
\end{nb}


\begin{defi}\label{tempsdevol}
For any $(x,v) \in \overline{\Omega} \times \mathbb{R}^N,$ define
$$t_{\pm}(x,v)=\sup\{\,t > 0\;;x \pm sv \in \Omega,\;\;\forall \,0< s <
t\,\}=\inf\{\,s > 0\,;\,x \pm sv \notin \Omega\}.$$ For the sake of
convenience, we will set
$$\tau(x,v):=t_-(x,v) + t_+(x,v) \:\:\text { for any } (x,v) \in \overline{\Omega} \times
\mathbb{R}^N.$$
\end{defi}
Let us define
$$\vartheta(x,v)=\int_{-t_+(x,v)}^{t_-(x,v)} \Sigma(x-sv,v) \d s,\qquad
\qquad (x,v) \in \overline{\Omega} \times \mathbb{R}^N.$$ One proves
easily the following thanks to Assumption \ref{asssig}.

\begin{lemme}\label{lemthe}For any $(x,v) \in \Omega \times \mathbb{R}^N$,
$$\vartheta(x,v)=\vartheta(x,-v)=\int_{0}^{\tau(x,v)} \Sigma(x+t_-(x,v)v-sv,v)\d s.$$
In particular, for any $(x,v) \in \Gamma_-$
$$\vartheta(x,v)=\int_0^{t_+(x,v)}\Sigma(x+sv,v)\d s.$$
\end{lemme}

\n Now let us investigate the resolvent of
$T$. For any $\lambda \in \mathbb{C}$ such that $\mathrm{Re} \lambda
> 0$, let us define $M_{\lambda} \in \mathscr{B}(L^p_-,L^p_+)$ by
\begin{equation*}
M_{\lambda}u(x,v)=u(x-\tau(x,v)v,v)\exp\left\{-\displaystyle
\int_0^{\tau(x,v)}\lambda+\Sigma(x-sv,v)\d s\right\},\qquad (x,v)
\in \Gamma_+,
\end{equation*}
and let $B_{\lambda} \in \mathscr{B}(L^p_-,X_p)$ be given by
\begin{equation*}
B_{\lambda}u(x,v)=u(x-t_-(x,v)v,v)\exp\left\{-\displaystyle
\int_0^{t_-(x,v)}\lambda+\Sigma(x-sv,v)\d s\right\},\quad (x,v) \in
\Omega \;.
\end{equation*}
In the same way, let $G_{\lambda} \in \mathscr{B}(X_p,L^p_+)$ be
given as
\begin{equation*}
G_{\lambda}\varphi(x,v)=\displaystyle\int_0^{\tau(x,v)}\varphi(x-sv,v)
\exp\left\{-\displaystyle
\int_0^{s}\lambda+\Sigma(x-tv,v)\d t\right\}\d s,\quad (x,v) \in
\Gamma_+\;;
\end{equation*}
and $C_{\lambda} \in \mathscr{B}(X_p)$ be defined as
\begin{equation*}
C_{\lambda}\varphi(x,v)=\displaystyle
\int_0^{t(x,v)}\varphi(x-tv,v)\exp\left\{\int_0^{t}\lambda+\Sigma(x-sv,v)\d
s\right\}\d t,\:\:\:(x,v) \in \Omega\;.
\end{equation*}

\n It is not difficult to show the following in the spirit of \cite{LM1}.

\begin{propo}\label{reso}
Let $0 < \gamma < 1$ be fixed and let $H \in
\mathscr{B}(L^p_+,L^p_-)$ be defined by:
\begin{equation}\label{eq-H}H(\phi_+)(x,v)=\gamma \phi_+(x,-v) \qquad \text{ for any }
\quad (x,v) \in \Gamma_-.\end{equation} If $\lambda \in \mathbb{C}$
is such that $1 \in \varrho(M_{\lambda}H)$, then $\lambda \in
\varrho(T)$ with \begin{equation}\label{resolvanteT}
(\lambda-T)^{-1}=B_{\lambda}H(I-M_{\lambda}H)^{-1}
G_{\lambda}+C_{\lambda}.\end{equation} In particular, if there is
$\lambda_0 \in \mathbb{R}$ such that
\begin{equation*}\label{resolvante}r_{\sigma}(M_{\lambda}H)
<1\qquad \qquad \forall \,\mathrm{Re} \lambda >
\lambda_0,\end{equation*} then $\left\{\lambda \in
\mathbb{C}\,;\,\mathrm{Re}\lambda
> \lambda_0\right\} \subset \varrho(T)$ and the resolvent of $T$ is
given by \eqref{resolvanteT}.
\end{propo}

\section{Study of the streaming operator and semigroup}

We shall focus in this section on the streaming operator associated
with bounce-back boundary conditions and the associated semigroup.

\subsection{Description of the spectrum of $T$}
To discuss the spectrum of $T$, we provide a more precise
description of the inverse operator of $(I-M_{\lambda}H) \in
\mathscr{B}(L^p_+)$. Precisely, let us define the measurable
function
$$
m_{\lambda}(x,v)=\gamma \,\exp\left\{-\int_{0}^{\tau(x,v)}
\lambda+\Sigma(x-sv,v)\d s\right\},\qquad (x,v) \in
\Gamma_+.\eqno $$

\n Before stating our first result we recall that the essential
range of the measurable function $m_{\lambda}(\cdot,\cdot)$, $\mathcal{R}_{\mathrm{ess}}(m_{\lambda})$,
is the set
$$\bigg\{ u\in \mathbb{C} \, :\, \big|\left\{ (x,v)\in \Gamma_{+};\,
\vert m_{\lambda}(x,v)-u\vert<\eps \right\}\big|\not= 0\ \ \forall
\eps>0\bigg\}$$

\n where $\vert A\vert$  denotes the Lebesgue measure of the set $A$.
Then, one has the following:
\begin{propo}\label{mlambdaH}
Let $\lambda \in \mathbb{C}$ be such that $1 \notin
\mathcal{R}_{\mathrm{ess}}(m_{\lambda})$. Then, $(I-M_\lambda H) \in
\mathscr{B}(L^p_+)$ is invertible with inverse given by
$$\left[(I-M_\lambda
H)^{-1}\psi\right](x,v)=(1-m_\lambda^2(x,v))^{-1}\left[(I+M_\lambda
H)\psi\right](x,v),\quad \forall\; \psi \in L^p_+,\quad (x,v) \in
\Gamma_+.$$
\end{propo}

\begin{proof} Let us fix $\lambda \in \mathbb{C}$  and consider the
equation: \begin{equation}\label{reso}\psi - (M_\lambda
H)\psi=g,\end{equation} where $g \in L^p(\Gamma_+)$ as well as the
unknown function $\psi.$ From \eqref{reso}, one sees that
$$\psi(x,v)-m_{\lambda}(x,v)\psi(x-\tau(x,v)v,-v)=g(x,v),\qquad
\qquad (x,v) \in \Gamma_+.$$ For any \textit{fixed} $(x,v) \in
\Gamma_+$, one has
\begin{multline*}
g(x-\tau(x,v)v,-v)=\psi(x-\tau(x,v)v,-v)-m_{\lambda}(x-\tau(x,v)v,-v) \times \\
\times \psi(x-\tau(x,v)v+\tau(x-\tau(x,v)v,-v)v,+v).\end{multline*}
Now, one sees easily that
$m_{\lambda}(x-\tau(x,v)v,-v)=m_{\lambda}(x,v)$ while
$$\psi(x-\tau(x,v)v+\tau(x-\tau(x,v)v,-v)v,+v)=\psi(x,v).$$
Therefore, one has
$$g(x-\tau(x,v)v,-v)=\psi(x-\tau(x,v)v,-v)-m_{\lambda}(x,v)\psi(x,v)$$
so that \begin{equation*}\begin{split} \left[M_{\lambda}H
g\right](x,v)&=m_{\lambda}(x,v)g(x-\tau(x,v)v,-v)\\
&=m_{\lambda}(x,v)\psi(x-\tau(x,v)v,-v)-m_{\lambda}^2(x,v)\psi(x,v).
\end{split}\end{equation*}
Since
$m_{\lambda}(x,v)\psi(x-\tau(x,v)v,-v)=\left[M_{\lambda}H\psi\right](x,v)$,
one obtains  from \eqref{reso} that
$$g(x,v)+\left[M_{\lambda}Hg\right](x,v)=
\left(1-m_{\lambda}^2(x,v)\right)\psi(x,v),\qquad (x,v) \in
\Gamma_+.$$ This leads to an explicit expression of the solution to
\eqref{reso}:
$$\psi(x,v)=(1-m_{\lambda}^2(x,v))^{-1}\left[\left(I+M_{\lambda}H\right)g\right](x,v).$$
Defining
$$\mathscr{R}_\lambda
g(x,v)=(1-m_{\lambda}^2(x,v))^{-1}\left[\left(I+M_{\lambda}H\right)g\right](x,v)$$
it is not difficult to see that $$1 \notin
\mathcal{R}_\mathrm{ess}(m_\lambda) \implies \mathscr{R}_\lambda \in
\mathscr{B}(L^p_+),$$ and the above calculations show that
$\mathscr{R}_\lambda=(I-M_{\lambda}H)^{-1}$.\end{proof}

The precise picture of the spectrum of $T$ is given by the
following, which is in the spirit of \cite{Lods2}
\begin{theo}
For any $k \in \mathbb{Z}$, let us define
$$F_k(x,v)=\dfrac{\log
\gamma-\vartheta(x,v)}{\tau(x,v)}-i\dfrac{2k\pi}{\tau(x,v)},\qquad
\quad \forall  (x,v) \in  {\Omega} \times \mathbb{R}^N.$$ Then,
$$\sigma(T)=\overline{\bigcup_{k \in \mathbb{Z}} \mathcal{R}_{\mathrm{ess}}(F_k)}$$
where $\mathcal{R}_{\mathrm{ess}}(F_k)$ stands for the essential
range of $F_k$.
\end{theo}

\begin{proof}  Let us begin with the   inclusion $\supset$. Given
$\lambda \in \mathcal{R}_{\mathrm{ess}}(F_k)$ $(k \in \mathbb{Z}).$
Let $\epsilon
> 0$ and define
$$\Lambda_{\epsilon}:=\{(x,v) \in \overline{\Omega} \times
\mathbb{R}^N\,;\,|\lambda-F_k(x,v)| \leq \epsilon\}.$$ By the
definition of $\mathcal{R}_{\mathrm{ess}}(F_k)$,
$|\Lambda_{\epsilon}| \neq 0$ for any $\epsilon > 0.$
For any integer $n \in \mathbb{N}$, define $B_n=\{(x,v) \in
 {\Omega} \times \mathbb{R}^N, \tau(x,v) \geq 1/n\}.$ Since
$\tau(x,v) \geq 0$ for a. e. $(x,v) \in  {\Omega} \times
\mathbb{R}^N,$ one has
$$\Lambda_{\epsilon}=\bigcup_{n} \bigg(B_n \cap \Lambda_{\epsilon}\bigg)$$
so that, for any $\epsilon > 0$, there exists $n(\epsilon) \in
\mathbb{N}$ such that $|B_{n(\epsilon)} \cap \Lambda_{\epsilon}|
\neq 0.$ Define
$$A_{\epsilon}:=B_{n(\epsilon)} \cap \Lambda_{\epsilon}, \qquad
\epsilon > 0$$ so that $|A_{\epsilon}|\neq 0$ for any $\epsilon
> 0,$ while
\begin{equation}\label{Esin}
\esin\big\{\tau(x,v)\;;\;(x,v) \in A_{\epsilon}\big\}\geqslant
1/n(\epsilon)> 0, \qquad \text{ for any } \epsilon >
0.\end{equation} Then, let $\chi_{\epsilon}$ stand for the
characteristic function of the measurable set $A_{\epsilon}$
$(\epsilon >0)$. One sees that
$$\chi_{\epsilon}(x,v)=\chi_{\epsilon}(x,-v), \text{ and } \quad
\chi_{\epsilon}(x+tv,v)=\chi_{\epsilon}(x,v)$$ for $t>0$ small
enough. Now, for any $\epsilon >0$, one can define
\begin{equation*}
\varphi_{\epsilon}(x,v)= \chi_{\epsilon}(x,v) \exp
\left\{-\int_0^{t_-(x,v)}\bigg(F_k(x,v)+\Sigma(x-sv,v)\bigg)\d
 s\right\}, \quad (x,v) \in \Omega \times \mathbb{R}^N.\end{equation*}
 One sees that $\varphi_{\epsilon} \in
L^{\infty}(\Omega \times \mathbb{R}^N)$ for any $\epsilon > 0,$
since

$$ \essu_{(x,v) \in A_{\epsilon}}
\left|\,\exp
\left\{-\int_0^{t_-(x,v)}\bigg(F_k(x,v)+\Sigma(x-sv,v)\bigg)\d
s\right\}\right| < \infty$$ by virtue of \eqref{Esin}. Moreover,
since $t_-(x,v)=0$ for any $(x,v) \in \Gamma_-$, one sees that
$${\varphi_{\epsilon}}_{|\Gamma_-} =\chi_{\epsilon}.$$ Given $(x,v)
\in \Gamma_-$ and $\epsilon >0$, one has
\begin{equation*}
\varphi_{\epsilon}(x,-v)= \chi_{\epsilon}(x,-v) \exp
\left\{-\int_0^{t_-(x,-v)}\bigg(F_k(x,-v)+\Sigma(x+sv,-v)\bigg)\d
s\right\}.\end{equation*} Since $t_-(x,-v)=t_+(x,v)=\tau(x,v)$  and
 $F_k(x,-v)=F_k(x,v)$ (see Lemma \ref{lemthe}), one has
\begin{equation*}
\varphi_{\epsilon}(x,-v)=\chi_{\epsilon}(x,v) \exp
\left\{-\tau(x,v)F_k(x,v)-\int_0^{t_+(x,v)}\Sigma(x+sv,-v)\,\d
s\right\}.\end{equation*} Now, one checks that
\begin{equation*}\begin{split}
-\tau(x,v)F_k(x,v)&=-\log \gamma + \vartheta(x,v) +2ik\pi\\
 &=-\log \gamma + 2ik\pi  + \int_0^{t_+(x,v)} \Sigma(x+sv,v)\d s \qquad
 \forall (x,v) \in \Gamma_-,\end{split}\end{equation*}
where we used again Lemma \ref{lemthe}.
Thus, one sees that, for any $(x,v) \in \Gamma_-$,

$$ \varphi_{\epsilon}(x,-v)=\dfrac{1}{\gamma}\,
\varphi_{\epsilon}(x,v)$$

\n which exactly means that $\varphi_{\epsilon}$ fulfils the boundary
 conditions \eqref{eq-H}. Finally, it is easy to see that $\varphi_{\epsilon} \in \D(T).$
 Define now the net $(\psi_\epsilon)_\epsilon$ by
 \begin{equation}\label{psi}
 \psi_{\epsilon}=\varphi_{\epsilon}/\|\varphi_{\epsilon}\|,
 \qquad \qquad \epsilon > 0\end{equation} and let $g_{\epsilon}=(\lambda-T)\psi_{\epsilon},$ i. e.
$$g_{\epsilon}(x,v)=\left(\lambda+\Sigma(x,v)\right)\psi_{\epsilon}(x,v)+
v \cdot \nabla_x \psi_{\epsilon}(x,v).$$
Using the fact, for any
$(x,v) \in \Omega \times \mathbb{R}^N,$ and for any $t>0$ small
enough, $\vartheta(x+tv,v)=\vartheta(x,v)$ (note that the same
occurs for $F_k$) one can prove that
$$g_{\epsilon}(x,v)=\big(\lambda-F_k(x,v)\big)\psi_{\epsilon}(x,v),
\qquad (x,v)\in \Omega \times \mathbb{R}^N, \quad \epsilon >0$$
and so
\begin{equation*}
\|g_{\epsilon}\| \leq \essu_{(x,v) \in A_{\epsilon}}
\left|\lambda-F_k(x,v)\right|\, \|\psi_{\epsilon}\| \leq \epsilon.
\end{equation*} This, together with \eqref{psi} achieves to show
that $\left(\psi_\epsilon\right)_{\epsilon >0}$ is a
\textit{singular net} of $T$ so that $\lambda \in \sigma(T).$ The
closedness of the spectrum ensures that
$$\overline{\bigcup_{k \in \mathbb{Z}} \mathcal{R}_{\mathrm{ess}}(F_k)} \subset
\sigma(T).$$

\n Let us prove now the converse inclusion. Assume that $\lambda \notin
\mathcal{R}_{\mathrm{ess}}(F_k)$ for any $k \in \mathbb{Z}.$ Then,
for any $k \in \mathbb{Z},$ there exists $\beta_k
>0$ such that
$$|\lambda -F_k(x,v)| \geq \beta_k \qquad \text{ a.e. } (x,v) \in \omm,$$
i.e.
$$\left|\dfrac{\log \gamma -\vartheta(x,v) - 2ik
\pi}{\tau(x,v)} -\lambda\right| \geq \beta_k \qquad \text{ a.e. }
(x,v) \in \omm.$$ Then,
$$|\log \gamma - \vartheta(x,v)-2ik\pi -\lambda\,\tau(x,v)|\geq
\tau(x,v) \,\beta_k
 \qquad \text{ a. e. } (x,v) \in \omm.$$ This means
that, for any integer $n \geq 0,$ there exists $c_n > 0$ such that
\begin{equation*}\label{cn1}
\left|\dfrac{\log \gamma -\vartheta(x,v) -
\lambda\,\tau(x,v)}{2\pi\,n} \pm i\right| \geq c_n \tau(x,v) \text{
a. e. } (x,v) \in \omm,\:n \geq 1,\end{equation*} and
\begin{equation*}\label{cn0}
|\log \gamma - \vartheta(x,v)-\lambda\,\tau(x,v)|\geq c_0\,\tau(x,v)
\qquad \text{ a. e. } (x,v) \in \omm.\end{equation*} Arguing as in
\cite{Lods2}, one can choose
$M>0$ such that
$$\left|\log \gamma- \vartheta(x,v)-\lambda \tau(x,v)\right|\leq M \qquad \text{
a. e. }  (x,v) \in \omm,$$ and
$$\left|\exp\left\{\log \gamma- \vartheta(x,v)-\lambda \tau(x,v)\right\}-1\right| \geq
\dfrac{C}{2}\prod_{n=1}^Nc_n^2\tau(x,v)^2 \qquad \text{ a. e. }
(x,v) \in \omm,$$ where
$$C=\esin_{(x,v) \in \omm}\left|\,\exp\left\{\dfrac{1}{2}(\log \gamma-
\vartheta(x,v)-\lambda \tau(x,v))\right\}\right|.$$ Moreover, one
can easily see that \begin{equation*}\label{est} \liminf_{\tau(x,v)
\to 0}|\exp\{\log \gamma- \vartheta(x,v)-\lambda \tau(x,v)\}-1| \geq
|1-\gamma| >0 .\end{equation*} In particular,

$$\label{1mlambda} \esin_{(x,v) \in \omm}|\exp\{\log \gamma-
\vartheta(x,v)-\lambda \tau(x,v)\}-1| >0.$$ This proves that $1
\notin \mathcal{R}_{\mathrm{ess}}(m_{\lambda})$. From Proposition
\ref{reso}, one gets that $(I-M_\lambda H)$ is invertible and
$\lambda \in \varrho(T)$.\end{proof}

\subsection{Explicit expression of the semigroup $\ut$}

We derive in this section the explicit expression of the semigroup
$\ut$ generated by the streaming operator $T$ associated to the
bounce-back boundary conditions $H$ given by \eqref{eq-H}. For
simplicity, we shall restrict ourselves to the case of a {\it
homogeneous collision frequency}:
$$\Sigma(x,v)=\Sigma(v),\qquad \qquad \forall (x,v) \in \Omega
\times \mathbb{R}^N.$$
\begin{theo}
The C$_0$-semigroup $\ut$ in $X_p$ generated by $T$ is given by
$$U(t) =\sum_{n=0}^\infty U_n(t), \qquad \forall t \geq 0,$$
where, for any fixed $t \geq 0$,
\begin{equation*}
\left[U_0(t)\varphi\right](x,v)= \varphi(x-tv,v)\exp(-\Sigma(v)
t)\chi_{\{t < t_-(x,v)\}},\qquad \varphi \in X_p, \quad (x,v) \in
\Omega \times \mathbb{R}^N\end{equation*} while, for any $n \geq 0$
$$\left[U_{2n+2}(t)\varphi\right](x,v)=\gamma^{2n+2}\exp(-\Sigma(v)
t)\chi_{\mathcal{I}_{2n+1}(x,v)}(t)\varphi\big(x-tv+(2n+2)\tau(x,v)v;v\big),$$
and
\begin{multline*}
\left[U_{2n+1}(t)\varphi\right](x,v)=\gamma^{2n+1}\exp(-\Sigma(v)
t)\chi_{\mathcal{I}_{2n}(x,v)}(t)\varphi\big(x+tv-2t_-(x,v)v- 2n
\tau(x,v)v,-v\big)\end{multline*} for any $\varphi \in X_p,$ and any
$(x,v) \in \Omega \times \mathbb{R}^N$, with
$$\mathcal{I}_k(x,v)=[k\tau(x,v)+t_-(x,v);(k+1)\tau(x,v)+t_-(x,v)],
\qquad \text{ for any  } k \in \mathbb{N}.$$
\end{theo}
\begin{proof} The proof is based upon the representation of the
resolvent \eqref{resolvanteT} and the use of the uniqueness of the
Laplace transform. Precisely, let $\lambda >0$ be fixed. According
to Proposition \ref{mlambdaH},
\begin{equation*}
\begin{split}
\left[(I-M_\lambda
H)^{-1}\psi\right](x,v)&=(1-m_\lambda^2(x,v))^{-1}\left[(I+M_\lambda
H)\psi\right](x,v)\\
&=\sum_{n=0}^\infty
\gamma^{2n}\exp(-2n(\lambda+\Sigma(v))\tau(x,v))\left[(I+M_\lambda
H)\psi\right](x,v),
\end{split}
\end{equation*}
for any  nonnegative $\psi \in L^p_+,$ i.e.
\begin{multline*}
\left[(I-M_\lambda H)^{-1}\psi\right](x,v) =\sum_{n=0}^\infty
\gamma^{2n}\exp(-2n(\lambda+\Sigma(v))\tau(x,v)) \psi(x,v) +\\
\sum_{n=0}^\infty
\gamma^{2n+1}\exp(-(2n+1)(\lambda+\Sigma(v))\tau(x,v))
\psi(x-\tau(x,v)v,-v),\qquad (x,v) \in \Gamma_+.\end{multline*} It
is then easy to see that, for any fixed $\varphi \in X_p$:
\begin{equation}\label{resoserie}B_{\lambda}H(I-M_{\lambda}H)^{-1}G_{\lambda}\varphi=\sum_{n=0}^\infty
\mathscr{J}_{ n }(\lambda)\varphi\end{equation}  where, for any $n
\geq 0$,
\begin{multline*}
\left[\mathscr{J}_{2n+1}(\lambda)
\varphi\right](x,v)=\gamma^{2n+1}\exp\left\{-2n(\lambda +
\Sigma(v))\tau(x,v)\right\}\exp\left\{-(\lambda+\Sigma(v))t_-(x,v)\right\}\times
\\\times \int_0^{\tau(x,v)}\varphi(x-t_-(x,v)v+sv,-v)\exp\left\{-(\lambda
+ \Sigma(v))s\right\} \d s,
\end{multline*}
and
\begin{multline*}
\mathscr{J}_{2n}(\lambda)
\varphi(x,v)=\gamma^{2n+2}\exp\left\{-(2n+1)(\lambda +
\Sigma(v))\tau(x,v)\right\}\exp\left\{-(\lambda+\Sigma(v))t_-(x,v)\right\}\times
\\\times \int_0^{\tau(x,v)}\varphi(x-t_-(x,v)v+\tau(x,v)v-s v, v)\exp\left\{-(\lambda
+ \Sigma(v))s\right\}   \d s.
\end{multline*}
 For fixed $(x,v) \in \Omega \times \mathbb{R}^N$, performing the change of
variable
$$t=2n\tau(x,v)+t_-(x,v)+s,\quad \d t=\d s,\quad t \in \mathcal{I}_{2n+1}(x,v)$$
in the above expression of $\mathscr{J}_{2n+1}(\lambda)$ leads
easily to
$$\left[\mathscr{J}_{2n+1}(\lambda)\varphi\right](x,v)=\int_0^\infty \exp(-\lambda t)
\left[U_{2n+1}(t)\varphi\right](x,v) \d t.$$ In the same way, the
change of variable
$$t=(2n+1)\tau(x,v)+t_-(x,v)+s,\quad \d t=\d s,\quad t \in \mathcal{I}_{2n+2}(x,v)$$
in the above expression of $\mathscr{J}_{2n+2}(\lambda)$ allows to
prove that
$$\left[\mathscr{J}_{2n+2}(\lambda)\varphi\right](x,v)
=\int_0^\infty \exp(-\lambda t)
\left[U_{2n+2}(t)\varphi\,\right](x,v)\, \d t, \qquad \forall n \geq
0.$$ Finally, it is easily seen that
$$\left[C_\lambda \varphi\right](x,v)=\int_0^\infty \exp(-\lambda t)
\left[U_0(t)\varphi\,\right](x,v)\, \d
t.$$ Therefore,
$$(\lambda-T)^{-1}\varphi=\sum_{n=0}^\infty \int_0^\infty \exp(-\lambda
t)U_n(t)\varphi\,\d t$$ for any $\varphi \in X_p$ for which the
series converges. Moreover, since $T$ generates a C$_0$-semigroup
$\ut$ in $X_p$, one also has
$$(\lambda-T)^{-1}\varphi=\int_0^\infty \exp(-\lambda
t)U (t)\varphi \, \d t.$$  From the uniqueness of the Laplace
transform, this yields
$$U(t)\varphi=\sum_{n=0}^\infty U_n(t)\varphi$$
for any nonnegative $\varphi \in X_p$ and, since all the operators
involved are clearly nonnegative and the positive cone of $X_p$ is
generating, the result holds for arbitrary $\varphi \in X_p$.
\end{proof}

\section{Spectral analysis of the perturbed semigroup}

We investigate now the spectral properties of the full semigroup
governing the problem \eqref{1}. Let us define the collision
operator $K$ by
$$K\varphi(x,v)=\int_{\mathbb{R}^N}\kappa(x,v,w)\varphi(x,w)\,\d w $$
where the kernel $\kappa(\cdot,\cdot,\cdot)$ is nonnegative over
$\Omega \times \mathbb{R}^N \times \mathbb{R}^N$. We shall assume
here that $K$ is a bounded operator, $K \in \mathscr{B}(X_p)$, $1
\leq p < \infty$, so that by the standard
 bounded perturbation theory, the operator $(T+K,\D(T))$ generates a
 C$_0$-semigroup $(V(t))_{t \geq 0}$ given by the following
 Dyson-Phillips expansion series:
$$V(t)=\sum_{j=0}^{\infty}V_j(t)$$
where $V_0(t)=U(t),$ $$V_j(t)=\displaystyle  \int_0^t
U(t-s)KV_{j-1}(s)\d s,\qquad (j \geqslant 1).$$

As indicated in the Introduction, our analysis does not cover the
case of transport equation in $L^1$-spaces, so we shall assume in
all this section that
$$1 < p < \infty.$$
Throughout the sequel, we shall assume that $K$ is a regular
operator in the following sense:

\begin{defi}\label{regular}
An operator $ {K} \in \mathscr{B}(X_p)$ $(1 < p < \infty)$ is said
to be regular if $ {K}$ can be approximated in the
 operator norm by operators of the form:
\begin{equation}\label{eq-4.1}
\varphi \in X_p \longmapsto \sum_{i \in
I}\alpha_i(x)\beta_i(v)\int_{\mathbb{R}^N} \theta_i(w)\varphi(x,w)\d
w \in X_p \end{equation}

\n where $I$ is finite, $\alpha_i \in L^{\infty}(\Omega)$, $\beta_i \in
L^p(\mathbb{R}^N,\d v)$ and $\theta_i \in L^q(\mathbb{R}^N,\d v)$,
$1/p+1/q=1.$\end{defi}

\begin{nb}\label{nbregular} Since $1 < p < \infty$, one notes
that the set $\mathcal{C}_c(\mathbb{R}^N)$
of continuous functions with compact support in $\mathbb{R}^N$ is
dense in $L^q(\mathbb{R}^N,\d v )$ as well as in
$L^p(\mathbb{R}^N,\d v)$ ($1/p+1/q=1$). Consequently, one may assume
in the above definition that $\beta_i(\cdot)$ and $\theta_i(\cdot)$
are continuous functions with compact supports in
$\mathbb{R}^N$.\end{nb}

\n We prove in this section the following compactness result,
generalizing known ones for $1D$-transport problems \cite{LS}

\begin{theo}\label{main} Assume $1 < p < \infty$. If $K \in \mathscr{B}(X_p)$
is a regular operator, then
the difference $V(t)-U(t)$ is compact for any $t \geq 0.$ As a
consequence, $\sigma_\mathrm{ess}(V(t))=\sigma_\mathrm{ess}(U(t))$
for any $t \geq 0$.
\end{theo}

\begin{nb} Notice that the compactness of the difference $V(t)-U(t)$
for any $t \geq 0$ implies that of
$(\lambda-T-K)^{-1}-(\lambda-T)^{-1}$
 for sufficiently large $\lambda.$ In particular,
$$\sigma_\mathrm{ess}(T+K)=\sigma_\mathrm{ess}(T).$$
This result was already obtained in \cite{Latrach} and  is valid for
more general  reentry boundary conditions.
\end{nb}

\n The rest of the paper is devoted to the proof of the above
Theorem. We shall adopt the so-called resolvent approach which
allows to infer the compactness of $$R_1(t)=V(t)-U(t), \qquad t \geq
0$$ from properties of the resolvent $(\lambda-T)^{-1}$ and $K$
only. The basis of our approach is a fundamental result owing  to M.
Sbihi \cite[Theorem 2.2, Corollary 2.1]{Sbihi} which, applied to our
case, asserts that, for $p=2$, if  $T$ is dissipative and there
exists $\alpha>w(U)$  ($w(U)$ denoting the type of the semigroup
$\ut$) such that
\begin{equation}\label{eq-4.2}\alpha +i\beta-T)^{-1}K(\alpha +i\beta-T)^{-1} \hskip 0.7 cm
\hbox{ is compact for all } \beta \in \R\end{equation} and
\begin{equation}\label{eq-4.3}\displaystyle\lim_{\beta\to \infty}(\Vert  K^*(\alpha
+i\beta-T)^{-1}K\Vert + \Vert  K(\alpha
+i\beta-T)^{-1}K^*\Vert)=0,\end{equation} then $R_1(t)=V(t)-U(t)$ is
compact on $X_2$ for all $t\geq 0.$

Notice that here, the streaming operator $T$ is dissipative on
$X_p$, $p\in (1,\infty)$, in particular for $p=2$. Moreover, the
compactness assumption \eqref{eq-4.2} follows  from Theorem 3.1 in
\cite{Latrach} and holds true for more general boundary conditions.
Therefore, we have only to check  that \eqref{eq-4.3} holds true
provided $K$ is a regular collision operator.\vs

\n Though M. Sbihi's result is a purely Hilbertian one, it has
already been noticed in \cite{LS} that it can be applied
successfully to neutron transport problems in $L^p$-spaces for any
$1 < p < \infty$. Actually,  since $K$ is regular and $R_1(t)$
depends continuously on $K \in \mathscr{B}(X_p)$, one may assume
that $K$ is of the form \eqref{eq-4.1} where, according to Remark
\ref{nbregular}, the functions $\beta_i$ and $\theta_i$ are
continuous with compact supports in $\mathbb{R}^N$. In this case,
$K$ is bounded in any $\mathscr{B}(X_r)$ and, by an interpolation
argument already used in \cite{LS},  we may restrict ourselves to
prove the compactness of $R_1(t)$ in $X_2$. Moreover, using a
domination argument as in \cite{LS}, there is no loss of generality
in proving the compactness of $R_1(t)$ in the special case
$$\Sigma(v)=\sigma >0, \qquad \gamma=1$$ Now, since $K$ is given by
\eqref{eq-4.1}, by linearity, Eq. \eqref{eq-4.3}, and consequently
Theorem \ref{main}, follow  from the following

\begin{lemme}\label{lem45} Let  $\beta_j, \theta_j$ be continuous functions
with compact support in $\mathbb{R}^N$ and  $\alpha_j \in L^\infty(\Omega)$,
 $j=1,2$. Then, there is some $\alpha > -\sigma$ such that

$$\displaystyle\lim_{|\beta| \to \infty}
 \|K_1(\alpha+i\beta-T)^{-1}K_2\|=0$$
where
 $$K_j\varphi(x,v)=\alpha_j(x)\beta_j(v)\int_{\mathbb{R}^N}
 \theta_j(w)\varphi(x,w)\d
 w,\qquad
 j=1,2; \quad \varphi \in X_2.$$\end{lemme}

\begin{proof} Let $\alpha >-\sigma$ be fixed. According to Eq.
\eqref{resolvanteT}
 $$(\lambda-T)^{-1}=B_{\lambda}H(I-M_{\lambda}H)^{-1}G_{\lambda}+C_{\lambda}.$$
Moreover, it is well-known from \cite{Mokhtar1} that
\begin{equation}\label{linearity}\displaystyle \lim_{|\beta| \to \infty}
 \|K_1C_{\alpha+i\beta} K_2\|=0, \qquad \forall \alpha > -\sigma.
 \end{equation} since $C_\lambda$ is the resolvent of the transport operator with
no-reentry boundary conditions. Therefore, one has to prove that
$\lim_{|\mathrm{Im}\lambda| \to \infty}
\|K_1B_{\lambda}H(I-M_{\lambda}H)^{-1}G_{\lambda} K_2\|=0$ where
$\mathrm{Re}\lambda >-\sigma.$ According to Eq. \eqref{resoserie},
it suffices to establish that
$$\lim_{|\beta| \to \infty}
 \|K_1\mathscr{J}_n(\alpha+i\beta) K_2\|=0, \qquad \forall n \in
 \mathbb{N}, \qquad \alpha > -\sigma.$$
Let us prove the result for $\mathscr{J}_{2n+1}(\alpha+i\beta)$, $n
\in \mathbb{N}$. Let $\lambda=\alpha+i\beta$, $\alpha
>-\sigma$  and let $n \in \mathbb{N}$ be fixed. Technical
calculations show that
$$K_1\mathscr{J}_{2n+1}(\lambda)K_2\varphi=\mathcal{A}_3
\mathcal{A}_2(\lambda)\mathcal{A}_1 \varphi$$ where
$$\mathcal{A}_1\::\:\varphi \in L^2(\Omega \times \mathbb{R}^N) \mapsto \mathcal{A}_1\varphi(x)=\alpha_2(x)\int_{\mathbb{R}^N}\varphi(x,w)
\theta_2(w)\d w \in L^2(\Omega),$$
$$\mathcal{A}_3\::\:\psi \in L^2(\Omega) \mapsto
\mathcal{A}_3\psi(x,v)=\alpha_1(x)\beta_1(v)\psi(x) \in L^2(\Omega
\times \mathbb{R}^N)$$ and $\mathcal{A}_2(\lambda)\::L^2(\Omega) \to
L^2(\Omega)$ is given by
\begin{multline*}
\mathcal{A}_2(\lambda)\varphi(x)=\int_{\mathbb{R}^N}
\exp\bigg\{-2n(\lambda+\sigma)
\tau(x,v')-(\lambda+\sigma) t_-(x,v')\bigg\}\beta_2(-v')\theta_1(v')\d v'\\
\int_0^{\tau(x,v')}\exp(-(\lambda+\sigma)
s)\varphi(x-t_-(x,v')v'+sv')\d s.
\end{multline*}
Therefore, it is sufficient to prove that
\begin{equation}\label{estim}
\lim_{|\beta| \to \infty}
 \|\mathcal{A}_2(\alpha+i\beta)
 \|_{\mathscr{B}(L^2(\Omega))}=0, \qquad \alpha >
 -\sigma.\end{equation}  To do so, we adopt the approach of \cite{Mokhtar1} and
\cite{Sbihi}. Precisely, setting $\mu=\lambda+\sigma$ and
$h(v)=\theta_1(v)\beta(-v),$ $v \in \mathbb{R}^N$, the change of
variable $s \mapsto t=s-\tau(x,v')$ leads to

$$\mathcal{A}_2(\lambda)\varphi(x)=\int_{\mathbb{R}^N}
h(v')\exp\{-2n\mu\tau(x,v')\}\d
v'\int_{-t_-(x,v')}^{t_+(x,v')}\varphi(x+tv')\exp(-\mu t) \d t.$$

\n Since $\Omega$ is convex, given $(x,v') \in \Omega \times
\mathbb{R}^N$,
$$t \in (-t_-(x,v'),t_+(x,v')) \Longleftrightarrow y=x+tv' \in
\Omega$$ The change of variable $y=x+tv'$ shows that
$\mathcal{A}_2(\lambda)$ is an integral operator :
$$\mathcal{A}_2(\lambda)\varphi(x)=\int_{\Omega}
\kappa(\lambda,x, y)\varphi(y)\d y$$ where
\begin{equation*}\begin{split}\kappa(\lambda,x,y)&=\int_{\mathbb{R}}
h\left(\frac{y-x}{t}\right)\exp\left\{-\mu t
-2n\mu\tau\left(x,\frac{y-x}{t}\right) \right\}\frac{\d t}{t^N}\\
&=\int_{\mathbb{R}} h\left(\frac{y-x}{t}\right)\exp\bigg\{-\mu t
-2n\mu \frac{|t|}{|x|}\tau\left(x+z,-\frac{x}{|x|}\right)
\bigg\}\frac{\d t}{t^N},\end{split}\end{equation*} where  we used
the know property $\tau(x,\frac{v}{s})=|s|\tau(x,v)$ for any $(x,v)
\in \Omega \times \mathbb{R}^N$ and any $s \in \mathbb{R}.$ Notice
that the very rough estimate
$$\|\mathcal{A}_2(\lambda) \|_{\mathscr{B}(L^2(\Omega))} \leq
\bigg(\int_{\Omega \times \Omega} |\kappa(\lambda,x,y)|^2 \d x \d
y\bigg)^{1/2}$$ apparently does not lead to \eqref{estim}. We have
to estimate the norm of $\mathcal{A}_2(\lambda)$ more carefully.
With respect to \cite{Mokhtar1}, one of the difficulty in estimating
$\|\mathcal{A}_2(\lambda)\|$ is that $\mathcal{A}_2(\lambda)$ is not
a convolution operator because of the dependence in $x$ of
$\tau(x,\cdot)$. To overcome this difficulty, we follow the approach
of \cite{Sbihi}. Precisely, set
$$N_\lambda(x,z)=\int_{\mathbb{R}}h\left(-\frac{x}{t}\right) \exp\left\{-\mu t -2n\mu\frac{|t|}{|x|}\tau\left(x+z,-\frac{x}{|x|}\right)\right\}
\frac{\d t}{t^N},$$ where $(x,z) \in \Omega \times \Omega$ with $x+z
\in \Omega$.  Let us point out that, from assumption, there is no
loss of generality assuming that there exist two constants $a, b
>0 $ such that
$$\mathrm{Supp} \left(h\right) \subset \left\{v \in \mathbb{R}^N\,;\,a \leq |v|
\leq b\right\}.$$

\n In this case, in the above integral, one sees that
$t \in \mathbb{R}$ is such that
$$a \leq \left|\frac{x}{t}\right| \leq b$$
which implies that $|t| \leq |x|/a $. This means that the above
integral over $\mathbb{R}$ reduces actually to an integral over
$\left[-\frac{d}{a},\frac{d}{a}\right]$ where $d$ is the diameter of
$\Omega.$ Then, $$\kappa(\lambda,x,y)=N_\lambda(x-y,y), \qquad
\text{ for any } (x,y) \in \Omega \times \Omega$$ and, setting
$$G_\lambda(x)=\sup_{z \in \Omega - x}
\left|N_\lambda(x,z)\right|,    \qquad x \in \Omega$$ one has
$$|\mathcal{A}_2(\lambda)\varphi(x)| \leq \int_{\Omega}
G_\lambda(x-y)|\varphi(y)|\d y, \qquad \forall x \in \Omega, \:
\varphi \in L^2(\Omega).$$ Consequently,
$$\left\|\mathcal{A}_2(\lambda)\right\|_{\mathcal{B}(L^2(\Omega))}
\leq \int_{\Omega}G_\lambda(x)\d x.$$

\n To prove \eqref{estim}, one has then to show that
$$\displaystyle \lim_{|\beta|
\to \infty} \int_{\Omega} G_{\alpha+i\beta}(x)\d x=0, \qquad \forall
\alpha > -\sigma.$$
 First, one sees that for any $(x,z)
\in \Omega \times \Omega$ with $x+z \in \Omega$ and any
$\lambda=\alpha+i\beta$, one has
$$\left|N_{\lambda}(x,z)\right| \leq \int_{\mathbb{R}}
\left|h\left(-\frac{x}{t}\right)\right|
\exp\big\{-(\alpha+\sigma)t\big\}\dfrac{dt}{t^N}$$ where we used the
fact that $\tau(\cdot,\cdot) \geq 0$. Then,
\begin{equation*}\begin{split}
\int_{\Omega} \sup_{\lambda=\alpha+i\beta}|G_\lambda(x)|\d x &\leq
\int_{\mathbb{R}}\exp\big\{-(\alpha+\sigma)t\big\}\dfrac{\d t}{t^N}\int_{\Omega}\left|h\left(-\frac{x}{t}\right)\right|\d x\\
& \leq
\int_{-\frac{d}{a}}^{\frac{d}{a}}\exp\big\{-(\alpha+\sigma)t\big\}\d
t \int_{\mathbb{R}^N} |h(v)|\d v < \infty,
\end{split}\end{equation*}
where we performed the change of variables $v=\frac{x}{t}$ in the
$x$ integral. Therefore, from the dominated convergence theorem, it
suffices to prove that

$$\lim_{|\beta| \to \infty} G_{\alpha+i\beta}(x)=0,\qquad \qquad
\text{ a. e. } x \in \Omega.$$

\n Using the fact that, for any fixed $x
\in \Omega$, the mapping $z \mapsto
\ds\tau\left(x+z,-\frac{x}{|x|}\right)$ is bounded, this can be done
as in \cite{Sbihi} thanks to the Riemann-Lebesgue's Lemma. This
achieves to prove  that
$$\lim_{|\beta| \to \infty}
 \|K_1\mathscr{J}_{2n+1}(\alpha+i\beta) K_2\|=0, \qquad \forall n \in
 \mathbb{N}, \qquad \alpha > -\sigma, \qquad \forall n \in \mathbb{N}.$$
One proves in the same way that
$$\lim_{|\beta| \to \infty}
 \|K_1\mathscr{J}_{2n}(\alpha+i\beta) K_2\|=0, \qquad \forall n \in
 \mathbb{N}, \qquad \alpha > -\sigma, \quad \forall n \in \mathbb{N}$$
 and, combined with \eqref{linearity} and \eqref{resoserie}, yield the result.
\end{proof}\vs

\begin{nb} Let us observe that Lemma \ref{lem45} allows  to describe  the asymptotic
spectrum of the transport operator $T+K$.  Indeed, combining  Lemma
\ref{lem45} with the compactness of  $(\lambda -T)^{-1}K$
\cite{Latrach}  and \cite[Lemma 1.1]{Mokhtar1} we infer that
\begin{enumerate}[i)\:]
\item  $\sigma (T+K)\cap \{\lambda \in \C \, :\, \mathrm{ Re
}\lambda>-\lambda^*\}$ consists of, at most, isolated eigenvalues
with finite algebraic multiplicity; \item for any $\eta>0$, the set
$\sigma (T+K)\cap \{\lambda \in \C \, :\, \mathrm{ Re
}\lambda>-\lambda^*+\eta\}$ is finite or empty.\end{enumerate}
Clearly, this result may be also derived from Theorem \ref{main}
using the fact that the spectral mapping theorem holds true for the
point spectrum.
\end{nb}

\end{document}